\documentclass[reqno,11pt,a4]{amsart}
\usepackage{amsmath,amssymb,tabularx,setspace,color,enumerate}
\usepackage{multirow}
\newtheorem{theorem}{Theorem}
\usepackage{bigstrut}
\usepackage{multirow}
\usepackage{booktabs,pdflscape,adjustbox}
\usepackage{hyperref}
\usepackage{array}
\usepackage{pdflscape}
\usepackage[margin=3.5cm]{geometry}
\usepackage{tabularx,graphics}

\newtheorem{definition}{Definition}
\newtheorem{remark}{Remark}
\newtheorem{corollary}{Corollary}
\makeatletter
\renewcommand\section{\@startsection {section}{1}{\z@}%
                                   {-3.5ex \@plus -1ex \@minus -.2ex}%
                                   {2.3ex \@plus.2ex}%
                                   {\normalfont\large\bfseries}}
\makeatother

\begin{document}


\title[]{\vspace{-0.2in}A new goodness of fit test for normal distribution based on Stein's characterization}
\author[]%
{   S\lowercase{udheesh} K K\lowercase{attumannil$^{\dag}$}\\
I\lowercase{ndian} S\lowercase{tatistical} I\lowercase{nstitute},  C\lowercase{hennai}, I\lowercase{ndia}.
}
\doublespace

    \begin{abstract} In this paper, we develop a simple non-parametric test for testing normal distribution based on the distance between empirical zero-bias transformation and empirical distribution.  The asymptotic properties of the test statistic are studied.  The finite sample performance of the proposed test is evaluated through a Monte Carlo simulation study.  The power of our test is compared with several other tests for normality.   We illustrate the test procedure using two real data sets.
\end{abstract}
\thanks{ {$^{\dag}$}{Author E-mail: \tt skkattu@isichennai.res.in}.}
\maketitle
\section{Introduction}
Test for normal distribution has great significance as most of the classical tests are developed on the assumption that the available data are generated from the normal distribution.  For the goodness of fit test associated with the normal distribution, we refer interested readers   to Thode  (2002), Schick et al. (2011), Shalit (2012),   Bera et al. (2016),  Torabi et al. (2006), Nikitin  (2018), Sulewski  (2019, 2020), Henze and  Jiménez-Gamero (2019), Betsch and Ebner (2020), Henze and Visagie (2020) and  Henze and Koch (2020) and the references therein. Among these, Betsch and Ebner (2020) developed a test for normal distribution based on the  distance between empirical zero-bias transformation and empirical distribution. Their test statistic has complicated expression and hence  the implementation of the test is difficult.    Motivated by Betsch and Ebner (2020) we develop a simple non-parametric test for testing normal distribution.

Stein's identity for normal distribution and its applications have been well studied in statistical literature. Let $X$ be a continuous random variable with finite mean $\mu$ and variance $\sigma^2$. Let $c(x)$ be a continuous function having first derivative. Then $X$ has normal distribution with mean $\mu$ and variance $\sigma^2$ if and only if
\begin{equation*}\label{stein}
E(c(X)(X-\mu))=\sigma^2E(c'(X)),
\end{equation*}provided the above expectations exist and the prime denotes the derivative with respect to $x$. Stein's type identity for general class of probability distributions and related  characterizations,  we refer interested readers  to Sudheesh (2009), Sudheesh and Tibiletti (2012) and Sudheesh and Dewan (2016) and the references therein.  Using Stein's characterization,  Betsch and Ebner (2021) developed a fixed point characterization for  normal random variables.  Using this fixed point characterization, we develop a  new goodness of fit test for  normal distribution. For recent devolvements on Stein's methods, we refer to  Anastasiou et al. (2023).

The rest of the article  is organized as follows.  In Section 2,  based on fixed point characterization, we develop a simple non-parametric test for testing normality. The  asymptotic properties of the test statistic are also studied.   A  Monte Carlo simulation study is carried out to evaluate  the finite sample performance of the proposed test and the results are reported in Section 3. We compare the power of our test with several  tests  for normality.  We illustrate our test procedure using two real data sets.  Concluding remarks  along with some open problems are given in Section 4.
\vspace{-0.2in}
\section{Test statistic}
\vspace{-0.1in} In this section, we develop a new goodness of fit for normal distribution. We use a fixed-point characterization for normal distribution to develop the test. Let $X$  be a continuous random variable with distribution function $F(.)$.  Assume  the mean  $\mu=E(X)$ and the variance $\sigma^2=Var(X)$  are finite. Define  $$e_X(x)=\frac{1}{\sigma^2}E(X(X-x)I(X\le x)),$$where $I$ denotes the indicator function.
We use the following characterization to develop the test for normal distribution.
\begin{theorem}
  (Betsch and Ebner, 2021). A continuous random variable $X$ with distribution function $F$,  $\mu=0$ and variance  $\sigma^2$ has the  normal distribution if  and only if $F(x)=e_X(x),\,\forall x\in R.$
  \end{theorem}\vspace{-0.2in}
Now  we consider the problem of  testing normality.
 Based on a random sample $X_{1}, ...,X_{n}$  from $F$, we are interested in  testing the null hypothesis
\begin{equation*}
  H_0:F\in \{N(0,\sigma^2);\sigma^2\in(0,\infty)\}
\end{equation*}\vspace{-0.1in}
against  a general alternatives
 \begin{equation*}\vspace{-0.1in}
   H_1:F\notin \{N(0,\sigma^2);\sigma^2\in(0,\infty)\}.
 \end{equation*}
\par For testing the above stated hypothesis, first we find a departure measure which discriminate the  null and alternative hypothesis. For this purpose, we consider  $ \Delta^*(F)$ given by
\begin{eqnarray*}\label{deltam}
  \Delta^*(F)&=&\int_{-\infty}^{\infty}(e_X(x)-F(x))dF(x).
  \end{eqnarray*}In view of Theorem 1, $\Delta^*(F)$ is zero under  $H_0$ and not equal to zero under $H_1$. Hence $\Delta^*(F)$  can be considered as a measure of departure  from  $H_0$ towards the alternative  hypothesis $H_1$.

We  propose  a test based on U-statistics. Hence, we write  $\Delta^*(F)$ as an expectation of the function of random variables. Consider
\begin{eqnarray}\label{del}
  \Delta^*(F)&=&\int_{-\infty}^{\infty}(e_X(x)-F(x))dF(x)\nonumber\\
&=&\int_{-\infty}^{\infty}(\frac{1}{\sigma^2}E(X(X-x)I(X\le x))-F(x))dF(x)\nonumber\\
&=&\frac{1}{\sigma^2}\int_{-\infty}^{\infty}\int_{-\infty}^{\infty}t(t-x)I(t\le x)dF(t))dF(x)-\frac{1}{2} \nonumber\\
&=&\frac{1}{\sigma^2}\int_{-\infty}^{\infty}\int_{-\infty}^{x}t(t-x)dF(t)dF(x)-\frac{1}{2}.
\end{eqnarray}
We observed that the probability density function of the random variable $\min(X_1,X_2)$ is $2\bar{F}(x)dF(x),$ where $\bar F(x)=1-F(x)$.   Hence  by  Fubini's theorem, we have
\begin{eqnarray}\label{eq2}
\int_{-\infty}^{\infty}\int_{-\infty}^{x}t^2dF(t)dF(x)&=& \int_{-\infty}^{\infty}t^2\int_{t}^{\infty}dF(x)dF(t)\nonumber\\
&=&\frac{1}{2}\int_{-\infty}^{\infty}2t^2\bar F(t)dF(t) \nonumber\\
&=&\frac{1}{2} E\big(\min(X_1,X_2)^2\big).
\end{eqnarray}\vspace{-0.1in}
Also
\begin{eqnarray}\label{eq3}
\int_{-\infty}^{\infty}\int_{-\infty}^{x}txdF(t)dF(x)&=&E\big(X_1X_2I(X_1<X_2)\big).
\end{eqnarray}
Substituting the equations (\ref{eq2}) and  (\ref{eq3}) in equation  (\ref{del}) we obtain
\begin{eqnarray}\label{deltaold}
  \Delta^*(F)=\frac{1}{2\sigma^2} E\big(\min(X_1,X_2)^2-2X_1X_2I(X_1<X_2)\big)-\frac{1}{2}.\nonumber
\end{eqnarray}Since the expected value of the random variables $\min(X_1,X_2)^2-2X_1X_2I(X_1<X_2)$ and $\min(X_1,X_2)^2-2X_1X_2I(X_1>X_2)$ are same, we consider the departure measure given by
\begin{eqnarray}\label{delta}
  \Delta(F)&=&\frac{1}{2\sigma^2} E\big(\min(X_1,X_2)^2-X_1X_2\big)-\frac{1}{2}.
\end{eqnarray}Next, we find an estimator of  the departure measure  given in (\ref{delta}). Let $h(X_1,X_2)=\min(X_1,X_2)^2- X_1 X_2$ be a symmetric kernel of degree 2. Then a U-statistic  given by
\begin{equation*}\label{USTAT}
  \widehat{\Delta}_1= \frac{1}{\binom{n}{2}} \sum_{i=1}^{n-1}\sum_{j=i+1}^{n}h(X_i,X_j),
\end{equation*}
is an unbiased and consistent estimator of  $\Delta_1(F)=E\big(\min(X_1,X_2)^2-X_1X_2\big)$. Also, the sample variance defined as $S^2=\frac{1}{n-1}\sum_{i=1}^{n}(X_i-\bar{X})^2$ is an unbiased and consistent estimator of $\sigma^2$.    Therefor, the  test statistic has the form
\begin{eqnarray*}
  \widehat\Delta&=&\frac{\widehat{\Delta}_1}{2S^2}-\frac{1}{2}.
\end{eqnarray*} Next, we express $ \widehat\Delta$ in a simple form. Let $X_{(i)}$, $i=1,\ldots,n$ be the $i$-th order statistic based on a random sample $X_1\ldots,X_n$ from $F$. Then   $\widehat\Delta$ can be written as
\begin{eqnarray}\label{deltahat}
  \widehat\Delta=\frac{1}{n(n-1)S^2} \sum_{i=1}^{n-1}(n-i)X_{(i)}^{2}- \frac{1}{n(n-1)S^2} \sum_{i=1}^{n-1}\sum_{j=i+1}^{n}X_{(i)}X_{(j)}-\frac{1}{2}.
\end{eqnarray}
The test procedure is to reject the null hypothesis  $H_0$ against the alternative hypothesis $H_1$ for large values of $\widehat\Delta$. We find a critical region of the test using  the simulated values of $\widehat\Delta$.
\begin{remark}Suppose $X_{1}, ...,X_{n}$ are random sample from $N(\mu,\sigma^2)$.  We can implement the proposed test based on the transformation $Y_i=(X_i-\bar{X})$, $i=1,\ldots,n$,  where $\bar{X}=\frac{1}{n}\sum_{i=1}^{n}X_i$ is the sample mean.
 \end{remark}
\noindent We  conduct a Monte Carlo simulation study to see whether the test is maintaining the size of the test while using the transformation $Y_i$. The result of simulation result is given in Table 1. From, Table 1, we can observe that the test has well-controlled type I error rates.
\begin{table}[h]
\begin{small}%
\caption{Empirical type I error}
\centering
\scalebox{0.99}{
\begin{tabular}{cccccccccc}
\hline
\multirow{2}{*}{} & \multicolumn{2}{c}{$N(1,2)$} & \multicolumn{2}{c}{$N(2,4))$} & \multicolumn{2}{c}{$N(3,5)$} & \multicolumn{2}{c}{$N(5,10)$} \\ \cline{1-9}
$n$& 1\% level& 5\% level& 1\% level& 5\% level& 1\% level& 5\% level&  1\% level& 5\% level   \\ \hline
10&	0.0108&	0.0502&	0.0111&	0.0508 & 0.0109&0.0512&0.0108&	0.0513\\
20&	0.0102&	0.0501&	0.0012&	0.0509&0.0105&	0.0508&	0.0106&	0.0512 \\
30&	0.0102&	0.0503&	0.0107&	0.0506&0.0103&	0.0506&	0.0102&	0.0508 \\
40&0.0100&	0.0501&	0.0102&	0.0501&0.0101&	0.0501&	0.0102&	0.0502 \\
50&	0.0101&	0.0501&	0.0101&	0.0500&0.0100&	0.0502&	0.0100&	0.0500 \\
\hline
\end{tabular}}\end{small}
\end{table}

\noindent Next, we obtain the asymptotic distribution of the test statistic.
\begin{theorem} As $ n \rightarrow \infty $,  $\sqrt{n}(\widehat{\Delta}_1-\Delta_1{(F)})$  converges in  distribution to a Gaussian random variable with mean zero and variance $4\sigma_{1}^{2}$, where $\sigma_{1}^2$ is given by
\begin{equation}\label{var}
\sigma_{1}^{2}=Var\Big(X^2\bar{F}(X)+\int_{-\infty}^{X}y^2dF(y)\Big).
\end{equation}
\end{theorem}
\noindent {\bf{Proof:}}  Using the central limit theorem on U-statistics,  we have the asymptotic normality of $\sqrt{n}(\widehat{\Delta}_1-\Delta_1(F))$ and the asymptotic  variance is  $4\sigma_{1}^{2}$, where $\sigma_{1}^{2}$ is given by (Lee, 1990)
\begin{equation}\label{var1}
\sigma_{1}^{2}=Var\big(E(h(X_{1},  X_{2})|X_1)\big).
\end{equation}Now, consider
\begin{eqnarray}\label{eq10}
    E(\min(X_1,X_2)^2|X_1=x) &=& E(x^2I(x<X_2)+X_2^2I(X_2\le x))\nonumber \\
   &=&   x^2\bar F(x)+\int_{-\infty}^{x}y^2dF(y).
\end{eqnarray}
Also
\begin{eqnarray}\label{eq11}
  E(X_1X_2|X_1=x)=
   0.
\end{eqnarray}
Substituting  (\ref{eq10}) and (\ref{eq11}) in (\ref{var1}) we obtain the asymptotic variance as specified in the equation (\ref{var}), which proves the theorem.\\
By applying Slutsky's theorem, from Theorem 2 we have following result.
\begin{corollary} As $ n \rightarrow \infty $,  $\sqrt{n}(\widehat{\Delta}-\Delta{(F)})$  converges in  distribution to  a Gaussian random variable with mean zero and variance $\sigma_{1}^{2}/\sigma^{4}$.
\end{corollary}\vspace{-0.1in}
\noindent Under the null hypothesis $H_0$ we know that $\Delta{(F)}=0$. Hence we have the following result. \vspace{-0.1in}
\begin{corollary}Under $H_0$, as $ n \rightarrow \infty $,  $\sqrt{n}\widehat{\Delta}$  converges in  distribution to  a Gaussian random variable with mean zero and variance $\sigma_{0}^{2}$, where $\sigma_0^2$ is given by
\begin{equation}\label{varnull}
\sigma_{0}^{2}=\frac{1}{\sigma^4}Var\Big(X^2\bar{F}(X)+\int_{-\infty}^{X}y^2dF(y)\Big).
\end{equation}
\end{corollary}
An asymptotic critical region of the test $\widehat{\Delta}$ can be constructed  using Corollary 2. Let $\widehat\sigma_{0}^2$ be a consistent estimator of $\sigma_{0}^2$. We reject the null hypothesis $H_{0}$ against the alternative hypothesis $H_{1}$ at a significance level $\alpha$, if
$\frac{ \sqrt{n} |\widehat{\Delta}| }{\widehat\sigma_0}>Z_{\alpha/2},$ where $Z_{\alpha}$ is the upper $\alpha$-percentile point of the standard normal distribution.

 Since  the distribution function  $F$ has no closed form for the normal distribution,  it is difficult to evaluate the null variance specified  in equation (\ref{varnull}). Hence we find the  critical region of the test based on Monte Carlo simulation. The critical points of the proposed test  for $\sigma^2=1$,  different sample sizes $n$ and  various significance levels $\alpha$ are given in Table 2. The critical  values are obtained based on  one million repetition. In  Table 2, lq and uq stands for lower and upper sample quantiles of the distribution of $\widehat\Delta$, respectively. We determine lower ($c_1$) and upper ($c_2$)  quantiles in such a way that $P(\widehat\Delta<c_1)=P(\widehat\Delta>c_2)=\alpha/2$.

\begin{table}[h]
\begin{small}
\caption{Simulated critical points of the test $\widehat{\Delta}$ for $N(0,1)$}
\centering\scalebox{0.9}{
\begin{tabular}{cccccccccc}
\hline
 &  \multicolumn{2}{c}{$\alpha =0.10$}&  \multicolumn{2}{c}{$\alpha =.05$}&  \multicolumn{2}{c}{$\alpha =0.02$}& \multicolumn{2}{c}{$\alpha =0.01$} \\ \hline
 $n$& lq  &uq&lq  &uq&lq  &uq&lq  &uq\\ \hline
 10&     -0.420&0.590&-0.472&0.752&-0.528&0.957&-0.566&1.106\\
 15&	 -0.357&0.466&-0.404&0.586&-0.456&0.735&-0.490&0.846\\
 20&	 -0.316&0.399&-0.360&0.498&-0.408&0.620&-0.439&0.712\\
 25&	 -0.286&0.352&-0.328&0.437&-0.373&0.541&-0.403&0.617 \\
 30&	 -0.264&0.319&-0.304&0.396&-0.346&0.489&-0.374&0.556 \\
 35&	 -0.245&0.293&-0.283&0.362&-0.325&0.445&-0.352&0.504\\
 40&	 -0.231&0.273&-0.267&0.336&-0.307&0.413&-0.333&0.467\\
 45&	 -0.220&0.256&-0.254&0.314&-0.292&0.387&-0.317&0.436 \\
 50&	 -0.209&0.242&-0.242&0.298&-0.278&0.364&-0.303&0.411 \\
 55&	 -0.200&0.230&-0.232&0.282&-0.268&0.344&-0.291&0.390\\
 60&	 -0.192&0.219&-0.223&0.269&-0.257&0.327&-0.280&0.369\\
 65&	-0.185 &0.211&-0.215&0.257&-0.248&0.327&-0.270&0.355 \\
 70&	 -0.179&0.202&-0.208&0.246&-0.240&0.300&-0.262&0.338 \\
 75&	 -0.173&0.195&-0.201&0.238&-0.234&0.290&-0.255&0.327\\
 80&	 -0.168&0.188&-0.196&0.230&-0.227&0.280&-0.247&0.314\\
 85&	 -0.163&0.183&-0.190&0.222&-0.221&0.270&-0.241&0.304 \\
 90&	 -0.159&0.176&-0.185&0.222&-0.215&0.261&-0.235&0.295 \\
 95&	 -0.155&0.172&-0.181&0.209&-0.210&0.261&-0.230&0.285\\
 100&	 -0.151&0.167&-0.177&0.204&-0.205&0.248&-0.224&0.277\\ \hline
\end{tabular}}\end{small}
\end{table}

\section{Empirical evidence and data analysis }
To evaluate the finite sample performance of the proposed  test, we conduct a Monte Carlo simulation study using R software. The simulation is repeated ten thousand times. To show the competitiveness of our test with the existing test procedures we  compare the empirical powers of the same. Finally, we illustrate our test procedure using two real data sets.

\subsection{Monte Carlo Simulation}
First, we find the empirical type I error of the proposed test. For finding the empirical type I error, we generated samples  from standard normal distribution and the result is given in  Table 3.
\begin{table}[h]
\caption{Empirical type I error of  $\widehat\Delta$}
\begin{tabular}{cccccc|c}\hline
$n$ &   1\% level  & 5\% level \\ \hline
10& 0.0101&0.0504\\
20&0.0100&0.0501\\
30&	0.0101& 0.0502 \\
40&0.0101 & 0.0500 \\
50&0.0100 & 0.0500 \\
 \hline
\end{tabular}
\end{table}

\begin{table}\begin{small}
\caption{Empirical Power of $\widehat\Delta$: Gumbel   distribution}
\centering\scalebox{0.99}{
\begin{tabular}{cccccccccc}
\hline
\multirow{2}{*}{} & \multicolumn{2}{c}{$\theta=0$} & \multicolumn{2}{c}{$\theta=0.2$} & \multicolumn{2}{c}{$\theta=0.4$} & \multicolumn{2}{c}{$\theta=0.6$} \\ \cline{1-9}
 $n$& 1\% level& 5\% level& 1\% level& 5\% level& 1\% level& 5\% level&  1\% level& 5\% level   \\ \hline
 10&0.0149& 0.1082&0.0256&02414
 &0.1082&0.4364&0.1642&0.6847\\
 20&0.0387&0.2601&0.3771&0.5993 &0.4332&0.8386& 0.7450 &0.9724\\
 30&0.0984&0.4173&0.4131&  0.7816& 0.7587 &0.9782&0.9663&1.0000\\
 40&0.1983& 0.5857 &0.5854& 0.9297&0.9510&0.9926&1.0000&1.0000\\
50& 0.3032&0.6844&0.7604&0.9708 &1.0000&1.0000&1.0000&1.0000\\ \hline
\end{tabular}}\end{small}
\end{table}

\begin{table}[h]\begin{small}%
\caption{Empirical Power of $\widehat\Delta$: Logistic distribution}
\centering
\scalebox{0.99}{
\begin{tabular}{cccccccccc}
\hline
\multirow{2}{*}{} & \multicolumn{2}{c}{$\theta=1$} & \multicolumn{2}{c}{$\theta=1.5$} & \multicolumn{2}{c}{$\theta=2$} & \multicolumn{2}{c}{$\theta=3$} \\ \cline{1-9}
$n$& 1\% level& 5\% level& 1\% level& 5\% level& 1\% level& 5\% level&  1\% level& 5\% level   \\ \hline
10&	0.2591&	0.5028&	0.3847&0.6132 & 0.5853&	0.7929&	0.8969&	0.9434\\
20&	0.3586&	0.5113&	0.4592&	0.6739&	0.7762&	0.8450&	0.9643&	0.9787\\
30&	0.3729&	0.5279&	0.5397&	0.6936&	0.8197&	0.8791&	0.9755&	0.9893\\
40&	0.4116&	0.5411&	0.5902&	0.7041&	0.8531&	0.9274&	0.9932&	0.9972\\
50&	0.4188&	0.5553&	0.5948&	0.7567&	0.8960&	0.9467&	0.9944&	0.9983\\
\hline
\end{tabular}}\end{small}
\end{table}

Next, by using empirical power we evaluate the finite sample performance of $\widehat\Delta$ against various alternatives. For finding empirical power, we simulate observations from  Gumbel ($\theta,\,1$), logistic ($\theta$), extreme value ($\theta$), and $t(n)$ distributions. The empirical power obtained for the above-mentioned alternatives is reported in Tables  4-7. From these tables, we observe that the empirical power of the test approaches one for large values of $n$. For Gumbel distribution, the empirical power is very high even for small sample sizes when the value of  $\theta$ is away from zero. From Table 7, we observe that the empirical power becomes very low as the degrees of freedom of $t$ distribution increase. This may be due to the fact that $t$ distribution becomes closer to the normal distribution as degrees of freedom increase.

\begin{table}[h]\begin{small}
\caption{Empirical Power of $\widehat\Delta$: Extreme value distribution}
\centering\scalebox{0.99}{
\begin{tabular}{cccccccccc}
\hline
\multirow{2}{*}{} & \multicolumn{2}{c}{$\theta=1$} & \multicolumn{2}{c}{$\theta=1.5$} & \multicolumn{2}{c}{$\theta=2$} & \multicolumn{2}{c}{$\theta=3$} \\ \cline{1-9}
$n$& 1\% level& 5\% level& 1\% level& 5\% level& 1\% level& 5\% level&  1\% level& 5\% level   \\ \hline
 10&0.1099	&0.3852 &0.1625&0.3972 &0.2226 &0.6318&0.8028&0.9043
\\
 20& 0.2473& 0.3852& 0.2313&0.5188&0.5218 &0.6982&0.9091&0.9468\\
 30&0.3168 & 0.4327&0.2770&0.6100&0.5863&0.7196&0.9245&0.9631\\
 40& 0.3323&0.4865  &0.3121&0.6769&0.6256&0.7559&0.9650&0.9735\\
50&0.3431& 0.4928&0.3608& 0.7317&0.6464&0.7971&0.9802&0.9920\\\hline
\end{tabular}}\end{small}
\end{table}

\begin{table}\begin{small}
\caption{Empirical Power of $\widehat\Delta$:  $t$ distribution}
\centering\scalebox{0.99}{
\begin{tabular}{cccccccccc}
\hline
\multirow{2}{*}{} & \multicolumn{2}{c}{$df=5$} & \multicolumn{2}{c}{$df=6$} & \multicolumn{2}{c}{$df=8$} & \multicolumn{2}{c}{$df=10$} \\ \cline{1-9}
 $n$& 1\% level& 5\% level& 1\% level& 5\% level& 1\%level& 5\% level&  1\% level& 5\% level   \\ \hline
 10&0.1335&0.2152     &0.1205&0.1928&0.0813 &0.1572&0.0512&0.1381\\
 20&0.1916 &0.2994 &  0.1661&0.2563 &0.1126 &0.1935&0.0850& 0.1452 \\
 30&0.2841&0.3657&0.1985 &0.3247&0.1295 &0.2389&0.0916 &0.1799\\
 40&0.3189&0.4331& 0.2132&0.3570&0.1431&0.2801&0.1049&0.2028\\
50&0.3577&0.4968 &0.2816&0.4113&0.1760&0.2857&0.1373&0.2215\\\hline
\end{tabular}}\end{small}
\end{table}

Next, we compare the power of the proposed test  with some  classical tests and recently developed tests  for normality. We consider  Anderson and Darling (AD), Kolmogorov–Smirnov (KS), Jarque-Berra (JB), Cramer-von Mises (CvM), Lilliefors (LF), Shapiro- Wilk (SW), Shapiro-Francia (SF) and SJ (Gel et al., 2007) tests for comparison.  We also consider recently developed tests like Hn (Torabi et al., 2016), MAD (Sulewski,
2019) and  $L(\alpha,\beta)$ (Sulewski, 2020) tests for comparison.  We briefly  discuss  these tests for the sake of completeness.

Let $\bar X=\frac{1}{n}\sum_{i=1}^{n}X_i$ and $S^2=\frac{1}{n-1}\sum_{k=1}^{n}(X_k-\bar{X})^2$ be the estimators of $\mu$ and  $\sigma^2$, respectively. Let $Z_i=\frac{X_i-\bar{X}}{S}$, $i=1,\ldots,n$ and  $h(x)=(x-1)^2/(x+1)^2$. We denote $Z_{(i)}$ as the ordered values of  $Z_i$.   Torabi et al. (2016) proposed a test statistic given by
$$H_n=\frac{1}{n}\sum_{i=1}^{n}h\left(\frac{1+\Phi(Z_{(i)})}{1+i/n}\right),$$
for testing normality, where $\Phi$ is the distribution function of the standard normal random variable.
Sulewski (2019) obtained a modified Anderson-Darling (MAD) test
for normality.  The MAD test has the form\vspace{-0.1in}
$$MAD=\sum_{i=1}^{n}\frac{(F_n(i)-\Phi(Z_{(i)}))^2}{\Phi(Z_{(i)})(1-\Phi(Z_{(i)}))}\phi(Z_{(i)}),$$
where $F_n(i)=(i-0.375)/(n+0.25)$ and $\phi$ is the derivative of $\Phi$.
Sulewski (2020) proposed a modified Lilliefors test. For $0\le \alpha,\beta\le 1$ consider the transformation (Blom, 1958)\vspace{-0.1in}
$$F_{\alpha,\beta}(X_{(i)})=\frac{i-\alpha}{n-\alpha-\beta+1}.$$
Using this transformation Sulewski (2020) obtained the $\alpha,\beta$ corrected  Lilliefors test statistics ($L(\alpha,\beta)$)  given by
$$L(\alpha,\beta)=\sup_{x}|F_{\alpha,\beta}(X_{(i)})-\Phi(Z_i)|.$$

   For power comparison, we consider a wide range of alternatives including standard Gumbel, standard logistic, $t$ and standard extreme value distributions. The forms of the density functions of the alternatives are given in Table 8. In Table 8, we use $\phi$ and $\Phi$ to denote standard normal density and distribution function, respectively. In power comparison, we also consider Tukey lambda distribution with quantile function given by $$Q(u,\lambda)=\begin{cases}
                                                                \frac{1}{\lambda}(u^{\lambda}-(1-u)^{\lambda}), & \mbox{if } \lambda\ne 0 \\
                                                                \log\left(\frac{u}{1-u}\right), & \mbox{if }   \lambda= 0
                                                              \end{cases}.$$
  We use R command for generating sample from all the alternatives  which are available in the R packages `rmutil', `extraDistr' and `sn'. In  simulation, we  consider the transformation $Y_i$ discussed in Remark 1 when the alternative is $N(5,245)$.
                                                              \begin{table}[h]
\caption{Alternative distributions used in the simulation study}
\centering
\begin{adjustbox}{width=\textwidth}
\begin{tabular}{lllll}
\hline
 Name of alternatives  &Probability density function &Notation \\\hline
  Extreme value&$\frac{1}{\beta}e^{\frac{x-\mu}{\beta}}\exp[-e^{\frac{x-\mu}{\beta}}]$&EV($\mu,\beta$)\\
 Gumbel  & $\frac{1}{\sigma}\exp\Big[\frac{\mu-x}{\sigma}-e^{(\frac{\mu-x}{\sigma})}\Big]$ &$GU(\mu,\sigma)$\\
Laplace &$(2b)^{-1}e^{-|x-a|/b}$&LAP($a,b$)\\
Logistic&${(1+e^{-(x-\theta)/\sigma})}{\sigma^{-1}(1+e^{-( x-\theta)/\sigma})^{-2}}$&LOS($\theta,\sigma$)\\
 Skew-normal &$\frac{2}{\sigma}\phi\big(\frac{\mu-x}{\sigma}\big)\Phi\big(c\frac{\mu-x}{\sigma}\big)$&$SK(\mu,\sigma,c)$\\
 Student's $t$&$(n\pi)^{-1/2}\Gamma^{-1}\left(\frac{n}{2}\right)\Gamma\left(\frac{n+1}{2}\right)\left(1+\frac{x^2}{n}\right)^{\frac{-(n+1)}{2}}$&$t(n)$\\
Uniform&$\frac{1}{b-a}I(a\le x\le b)$&U$(a,b)$\\
Compound Gumbel &$pGU(\mu_1,\sigma_1)+(1-p)GU(\mu_2,\sigma_2)$&CGU($\mu_1,\sigma_1,\mu_2,\sigma_2,p$)\\
 Compound normal &$pN(\mu_1,\sigma_1)+(1-p)N(\mu_2,\sigma_2)$&CND($\mu_1,\sigma_1,\mu_2,\sigma_2,p$)\\
  Normal Gumbel&$pN(\mu_1,\sigma_1)+(1-p)GU(\mu_2,\sigma_2)$&NGUM($\mu_1,\sigma_1,\mu_2,\sigma_2,p$)\\
\hline
 \end{tabular}\end{adjustbox}
\end{table}
\vspace{-0.05in}

 The similarity measure defined by Sulewski  (2020) can be used to  compare the standard normal distribution to an alternative distribution. Suppose that the alternative distribution has probability density function $f^*(x)$. Then the similarity measure is given by
\begin{equation*}
  M=\int_{-\infty}^{-\infty}\min\left\{\phi(x),f^*(x)\right\}dx.
\end{equation*}

\noindent Note that $0\le M\le 1$ and $M$ has value one  when both the densities are identical.     We illustrate it for a particular alternative  considered for simulation. When the alternatives are $t$ distributions with  5, 8 and 10 degrees of freedom, we obtain the values of $M$  as 0.939, 0.961 and 0.969, respectively.

The results of the Monte Carlo simulation study are reported  in Table 9. From Table 9, we observe that the size of all the tests we consider are  approximately equal to 0.05, hence it is comparable. The power of all the test increases when sample size increases.   Except few alternatives, our test has better power compared to the other tests.  For $U(-1.5,1.5)$ our test has poor performance, however in the case of  $U(-1,1)$ our test gives very good power for $n=50$. For $t$ distribution power of the proposed  test is better than that of the other tests. We also note that the power of all tests decreases as the degrees of freedom  of the $t$ distribution increases. For skew normal distribution our test has power  one even for $n=25$. In the case of compound Gumbel distribution  the power is very high compared to other tests. Overall, the  results given in Table 9 support  the claim that the proposed test is very  good alternative for several classical tests and recently obtained tests for normality.

 \begin{landscape}
\begin{table}\begin{small}
\caption{Power comparison:  Different alternatives ($\alpha=0.05$)}
\centering\scalebox{0.92}{
\begin{tabular}{cccccccccccccccccccccccccccc}
\hline
\multirow{2}{*}{} & \multicolumn{2}{c}{} &  \multicolumn{2}{c}{} & \multicolumn{2}{c}{} & \multicolumn{2}{c}{}& \multicolumn{1}{c}{} & \multicolumn{2}{c}{} & \multicolumn{2}{c}{$L(\alpha,\,\beta)$} &  \multicolumn{2}{c}{} & \\ \cline{12-16}
{Alternatives}& $n$&  AD& CVM& SW& SF& SJ&$\chi^2$& JB & Hn & MAD&(0,1)&(1,0)&(0,0)&(0.5,0.5)&(1,1)&$\widehat\Delta$ \\ \hline
 \multirow{2}{*}{N(0,1)}&25& 0.0531 &0.0521 &0.0532& 0.0532 &0.0510&0.0515&0.0514&0.0510& 0.0512&0.0561 &0.0496& 0.0512& 0.0521&  0.0512 &0.0499 \\
                        &50 &0.0498 &0.0489& 0.0500& 0.0501 &0.0508&0.0510&0.0492&0.0500&0.0505 &0.0479& 0.0503& 0.0488 &0.0489 &0.0490 &0.0502\\
\multirow{2}{*}{N(5,245)} &25&0.0406& 0.0406 &0.0450& 0.0450& 0.0500& 0.0505& 0.0409 &0.0408& 0.0490& 0.0500& 0.0459 &0.0459 &0.0459 &0.0409 &0.0489\\
 &50&0.0550& 0.0502& 0.0502& 0.0502& 0.0550 &0.0452& 0.0501& 0.0459& 0.0470& 0.0551 &0.0551& 0.0501
& 0.0459& 0.0490& 0.0496
\\
\multirow{2}{*}{EV(0,1)} &25& 0.3544& 0.3191& 0.3958& 0.3958& 0.2592& 0.3718& 0.3199& 0.2207& 0.1536& 0.1263& 0.3302& 0.2448& 0.2571& 0.2531& 0.7920\\
&50 &0.6047& 0.5469& 0.6783& 0.6783& 0.3437& 0.6284& 0.5462& 0.4450& 0.3961& 0.3049& 0.5331& 0.4542 & 0.4630& 0.4697& 0.9585\\
\multirow{2}{*}{SK(0,1,5)} &25 &   0.3163& 0.2843& 0.3574& 0.3574& 0.1381& 0.2593& 0.2843& 0.3934& 0.5886& 0.3103 & 0.1001& 0.2362& 0.2312& 0.2122& 1.0000\\
 &50&0.5806& 0.5315& 0.6597& 0.6597& 0.1592& 0.4805& 0.5315& 0.6166& 0.7968& 0.5055&  0.2773& 0.4414& 0.4274& 0.4174&1.0000\\
 \multirow{2}{*}{GU(0,1)} &25 &0.3383 &0.3022 &0.3864& 0.3864 &0.2409& 0.3326 &0.3020 &0.3883 &0.1568& 0.3321 &0.1309& 0.2624 & 0.2580 &0.2585& 0.3880\\

                          &50 &0.6125 &0.5659& 0.7139 &0.7139& 0.3630 & 0.6511& 0.5659 &0.6348& 0.4016& 0.5456& 0.2981& 0.4381& 0.4401 &0.4523&  0.7147\\
 \multirow{2}{*}{t(5)}
& 25& 0.2066& 0.1792 &0.2238& 0.2238 &0.2917& 0.2727& 0.1791& 0.1788& 0.0799& 0.1433& 0.1461 &0.1288& 0.1452& 0.1746 &0.3468\\
&50& 0.3171& 0.2807 &0.3782 &0.3782& 0.4555 &0.4304& 0.2804 &0.2682& 0.1447& 0.2114& 0.2097 &0.1892 & 0.2114 &0.2314 &0.4922\\
 \multirow{2}{*}{t(8)} &25&0.1236& 0.1098& 0.1324& 0.1324& 0.2038& 0.1743& 0.1092& 0.1016 &0.0603 &0.0929 &0.0825& 0.0858 & 0.0988& 0.1112& 0.2226\\
&50&0.1713& 0.1483& 0.2166& 0.2166& 0.2771& 0.2749& 0.1484& 0.1387& 0.0883 &0.1182& 0.1038& 0.0981& 0.1140& 0.1323& 0.2820\\
 \multirow{2}{*}{t(10)} &25& 0.1023& 0.0890& 0.1026 &0.1026& 0.1608& 0.1302& 0.0893& 0.0888& 0.0599& 0.0722& 0.0628& 0.0709& 0.0729& 0.0848& 0.1584\\
 &50&0.1138& 0.1090& 0.1443& 0.1443& 0.2149& 0.1950& 0.1092& 0.1090& 0.0582& 0.0857& 0.0839& 0.0724& 0.0809& 0.0939& 0.2260\\
 \multirow{2}{*}{U(-1,1)}
&25&0.2149&	0.1771&	0.2570&	0.2570&	0.0121&	0.0129&	0.1772&	0.1844&	0.5673&	0.1120&	0.1092&	0.1437&	0.1026&	0.0789&	0.1346\\
&50&0.5917&	0.4562&	0.7571&	0.7571&	0.0224& 0.0218&	0.4562&	0.4818&	0.8649&	0.2667&	0.2515&	0.3358&	0.2731&	0.2141&	0.9933\\
\multirow{2}{*}{U(-1.5,1.5)}
&25&0.2596&	0.2031&	0.3080&	0.3080&	0.0002&	0.0014&	0.2038&	0.2225&	0.2850&	0.1353&	0.1252&	0.1690&	0.1338&	0.0914&	0.0041\\
&50&0.5968&	0.4675&	0.7722&	0.7722&	0.0004&	0.0076&	0.4676&	0.4939&	0.6192&	0.2996&	0.2629&	0.3547&	0.2894&	0.2228&	0.0129\\
 \multirow{2}{*}{LOS(0,1)} &25&0.1183&	0.1130&	0.1282& 0.1282&	0.1837& 0.1671&	0.1130&	0.1130&	0.0079&	0.0896&	0.0918&	0.0812&	0.0896&	0.1027& 0.8462 \\
&50&0.1746&	0.1623&	0.1865&	0.1866&	0.2688&	0.2324&	0.1624&	0.1589&	0.0149&	0.1222&	0.1344&	0.1246&	0.1330&	0.1462&	0.9784\\
 \multirow{2}{*}{LAP(1,2)} &25&0.3029& 0.2892& 0.3003& 0.3003& 0.4815& 0.3534& 0.2893& 0.2712& 0.1190& 0.2377& 0.2173& 0.1949& 0.2372& 0.2750& 0.4658\\
 &50&0.5464& 0.5449& 0.5307& 0.5307& 0.7383& 0.5508& 0.5449& 0.5162& 0.2808& 0.4228& 0.4349& 0.3961&
 0.4427& 0.4741& 0.6720\\
 \multirow{2}{*}{TLD(1)} &25&0.2176& 0.1729& 0.2636& 0.2636& 0.0012&  0.0016& 0.1729& 0.1856& 0.5840& 0.1185& 0.1191& 0.1596& 0.1203& 0.0852& 0.1172\\
 &50&0.5857& 0.4593& 0.7558& 0.7558& 0.0022& 0.0087&  0.4593& 0.4912& 0.8837& 0.2470& 0.2732& 0.3241& 0.2718& 0.2194& 0.9927\\
 \multirow{2}{*}{NGUM(0,1,0,1,0.2)} &25& 0.2677& 0.2428&0.3057& 0.3057& 0.1818 & 0.2657& 0.2428 &0.3177 &0.1119& 0.2587& 0.0729& 0.1798& 0.1868& 0.1848& 0.2627\\
&50 &0.4885& 0.4515 &0.5514& 0.5514& 0.2707& 0.4945 &0.4515& 0.4985 &0.3017& 0.4116& 0.2458 &0.3467& 0.3546& 0.3516 &0.5554\\
 \multirow{2}{*}{CND(0,1,-4,1,0.7)} &25
&0.0480 &0.0501 &0.0430& 0.0430& 0.0621& 0.0450& 0.0501& 0.0571& 0.0480& 0.0561&
0.0470 &0.0531& 0.0521& 0.0501& 0.2442\\
&50& 0.0347& 0.0421& 0.0368& 0.0368& 0.0568& 0.0484& 0.0421 &0.0347& 0.0453&
0.0432& 0.0484& 0.0432& 0.0474& 0.0463& 0.2547\\
\multirow{2}{*}{CGU(0,1,2,1,0.8)} &25&0.3060&  0.2712& 0.3638&  0.3638& 0.2499& 0.3346& 0.2712& 0.3583& 0.1382& 0.3049& 0.1079& 0.2107& 0.2169& 0.2191& 0.4313\\
&50&0.5962& 0.5389& 0.6752 & 0.6752& 0.3351& 0.6030& 0.5385& 0.6178& 0.3664& 0.5071& 0.2826& 0.4219& 0.4277& 0.4239& 0.7408\\
 \multirow{2}{*}{CGU(0.8,1,2.7,1,0.8)} &25&0.3357& 0.3037& 0.3966& 0.3966 &0.2308 &0.3367& 0.3037 &0.3856& 0.1518& 0.3267& 0.1249& 0.2338& 0.2557& 0.2647& 0.9990\\
 &50&0.6016 &0.5596& 0.6957 &0.6957& 0.3514& 0.6146& 0.5596 &0.6336& 0.3944& 0.5355& 0.2843 &0.4334& 0.4364& 0.4384& 1.0000\\\\
 \hline
 \end{tabular}}\end{small}
\end{table}
\end{landscape}
\subsection{Data Analysis}
We illustrate  the implementation  of our test procedure using two real data sets. First we consider `Oxboys' data  from the R package `nlme'. This is a data on the height of a selection of 234 boys from Oxford, England versus a standardized age. We  consider height in cm to test the normality. The graphical visualization of the distribution is shown using histogram and is given in Figure 1.  The calculated value of the test statistic $\widehat\Delta$  is   -0.0126, which suggest that the underlying data follows normal distribution and the result is same with  the results obtained by many others.  We also  find the p-values for different tests and is reported  in Table 10.  From Table 10, we observe that  the p-value corresponds to $\widehat{\Delta}$, MAD and $\chi^2$ tests are higher compared to the other tests.

Next, we consider a data set named `Frets' presented in the R package `boot'. The data contains the measurements of the head length and breadth in millimeters of pairs of adult brothers in 25 randomly sampled families. The head length of the eldest son is considered for the test for normality. The distribution of the data is visualized using histogram (Figure 2), suggest the normality. The calculated value of the test statistic $\widehat\Delta$   is  0.0199, which confirms that the underlying data follows normal distribution.
The p-values for all the tests are given in Table 11. From Table 11, we observe that the $\widehat{\Delta}$, MAD and JB tests have higher p-values compare to other tests.
\begin{figure}[!htb]
   \begin{minipage}{0.48\textwidth}
     \centering
     \includegraphics[width=1.1\linewidth]{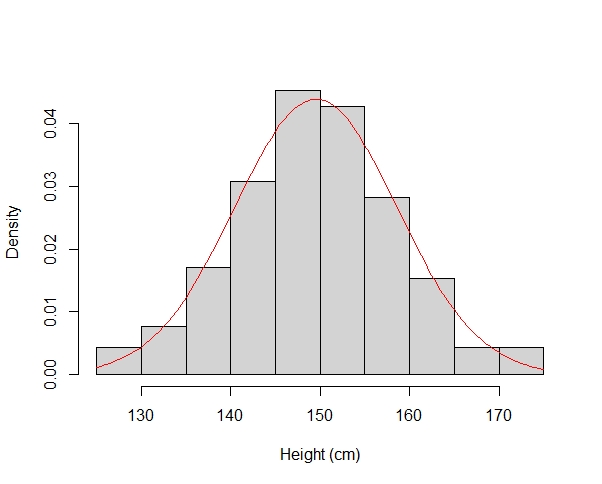}
     \caption{Histogram of the height of boys in Oxford, England}\label{Fig:Data1}
   \end{minipage}\hfill
   \begin{minipage}{0.48\textwidth}
     \centering
     \includegraphics[width=1.1\linewidth]{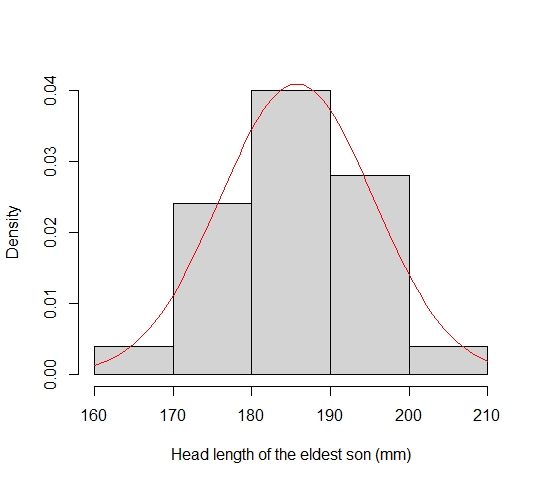}
     \caption{Histogram of the head length of the eldest son }\label{Fig:Data2}
   \end{minipage}
\end{figure}
\vspace{0.05in}
\begin{table}[ht]
  \centering
  \caption{P-values: Example I}
 \begin{adjustbox}{width=\textwidth}
  \begin{tabular}{ccccccccccccccccccccc}
    \hline
 Name of the test&   AD& CVM& SW& SF&$\chi^2$& JB & Hn  \\\hline
P-value& 0.7076& 0.6659& 0.6712& 0.7309& 0.9355& 0.7535 &0.6705 \\\vspace{-0.09in}
&&&&&&&&\\
Name of the test& MAD&LF(0,1)&LF(1,0)&LF(0,0)&LF(0.5,0.5)&LF(1,1)&$\widehat\Delta$\\\hline
P-value&1.0000&0.3109& 0.4723& 0.4075& 0.3835 &0.3636& 0.7760\\
    \hline
  \end{tabular}\end{adjustbox}
\end{table}
\begin{table}[ht]
  \centering
  \caption{P-values: Example II}
 \begin{adjustbox}{width=\textwidth}
  \begin{tabular}{ccccccccccccccccccccc}
    \hline
  Name of the test
&  AD& CVM& SW& SF& $\chi^2$& JB & Hn \\\hline
P-value&0.6434& 0.6800& 0.8407& 0.5637 &0.5323 &0.9785& 0.7579\\\vspace{-0.09in}
&&&&&&\\
  Name of the test&MAD&LF(0,1)&LF(1,0)&LF(0,0)&LF(0.5,0.5)&LF(1,1)&$\widehat\Delta$\\ \hline
 P-value&1.0000& 0.7846& 0.3736& 0.5463 &0.5612 &0.6050& 0.8885\\
    \hline
  \end{tabular}\end{adjustbox}
\end{table}
\noindent The R program for calculating the test statistic, the critical value and the corresponding p-value is given in  Appendix A.
\section{Concluding Remarks}
Stein's identity for normal random variables and its applications in different areas have been explored by many researchers. Based on fixed point characterization derived from Stein's identity, we developed a goodness of fit test for normal distribution. The test statistic has an asymptotic normal distribution.   Through a Monte Carlo simulation study, we have shown that the proposed test has good power for various alternatives.  The simulation study also shows that the test is compared with some well-known classical tests available in the literature. Finally, we illustrated our test procedure using two real data sets.

  Recently, based on Stein's method, Betsch and  Ebner (2021) obtained fixed point characterizations for a large class of absolutely continuous probability distributions. Using these fixed-point characterization, we can develop goodness of fit tests for different univariate continuous distributions. This work and recent papers by  Betsch and  Ebner (2020, 2021) are good starting points in this direction.


\vspace{-0.3in}
\section*{Appendix A}
\vspace{-0.2in}
\#\#..........The R program for finding test statistics, critical point and p-value.......\\
alfa=0.05 \quad\quad       \#significance level\\
delta=function(x)\{\\
s=0; n=length(x);\\
for(i in 1:(n-1))\{\\
for(j in (i+1):n)\{\\
s=s+(1/(2*var(x)))*((min(x[i],x[j])$)^2$-x[i]*x[j])\}\}\\
return((s/choose(n,2))-0.5)\}\\
library(nlme)\\
data=Oxboys\$height  \quad  \# input the data\\
data=data-mean(data)\\
delta(data)      \quad  \# statistic\\
n=length(data)\\
B=1e6    \,\,\,\,\,         \# number of experiments\\
z=matrix(rnorm(n*B,0,sqrt(var(data))),nrow=B, ncol=n)\\
H0=apply(z, 1, delta)\\
Q1=quantile(H0, alfa/2); Q1 \,\, \,\,\,\# critical value\\
Q2=quantile(H0, 1-alfa/2); Q2 \,\,\# critical value\\
sum(abs(H0)$>$abs(delta(data)))/B \,\,\, \,\#  p-value\\
\#\#.....................................................................

\vspace{-0.1in}
\begin{thebibliography}{xx}\vspace{-0.1in}


\bibitem{} Anastasiou, A., Barp, A., Briol, F. X., Ebner, B., Gaunt, R. E., Ghaderinezhad, F., ... and Swan, Y. (2023). Stein's method meets computationa statistics: a review of some recent developments. {\em Statistical Science},  38, 120--139.

\bibitem{} Bera, A. K., Galvao, A. F., Wang, L. and Xiao, Z. (2016). A new characterization of the normal distribution and test for normality. {\em Econometric Theory},  32, 1216--1252.



\bibitem{}  Betsch, S.  and Ebner, B. (2020). Testing normality via a distributional fixed point property in the Stein characterization. {\em TEST}, 29, 105--138.
\bibitem{} Betsch, S. and  Ebner, B. (2021). Fixed point characterizations of continuous univariate probability distributions and their applications. {\em Annals of the Institute of Statistical Mathematics}, 73, 31--59.
    \bibitem{blom} Blom, G. (1958). {\em Statistical Estimates and Transformed Beta Variables}. Wiley, New York.

\bibitem{} Henze, N.,  Jiménez-Gamero, M. D. (2019). A new class of tests for multinormality with iid and garch data based on the empirical moment generating function. {\em TEST}, 28, 499--521.

\bibitem{HK} Henze, N. and Koch, S. (2020). On a test of normality based on the empirical moment generating function. {\em Statistical Papers}, 61, 17--29.
\bibitem{} Henze, N. and Visagie, J. (2020). Testing for normality in any dimension based on a partial differential equation involving the moment generating function. {\em Annals of the Institute of Statistical Mathematics}, 72, 1109--1136.
\bibitem{} 	
Jing, B. Y., Yuan, J. and Zhou, W. (2009). Jackknife empirical likelihood. {\em Journal of the American Statistical Association}, 104, 1224--1232.


\bibitem{} Lee, A. J. (1990).  \textit{U-Statistics: Theory and Practice}, Marcel Dekker Inc., New York.


\bibitem{}Nikitin, Y. Y. (2018). Local exact Bahadur efficiencies of two scale-free tests of normality based on a recent characterization. {\em Metrika,} 81, 609--618.
    \bibitem{} Shalit, H. (2012). Using OLS to test for normality. {\em  Statistics \& Probability Letters,} 82, 2050--2058.

\bibitem{} Schick, A., Wang, Y. and Wefelmeyer, W. (2011). Tests for normality based on density estimators of convolutions. {\em Statistics \& probability letters}, 81, 337--343.

\bibitem{} Sudheesh, K. K. (2009). On Stein's identity and its applications. {\em Statistics \& Probability Letters,} 79, 1444--1449.
\bibitem{} Sudheesh, K. K. and Tibiletti, L. (2012). Moment identity for discrete random variable and its applications. {\em Statistics}, 46, 767--775.

\bibitem{}Sudheesh, K. K. and  Dewan, I.  (2016). On generalized moment identity and its applications: A unified approach. {\em Statistics,} 50, 1149--1160.

\bibitem{}Sulewski, P. (2019). Modification of Anderson-Darling goodness-of-fit test for normality. {\em Afinidad}, 76, 270--277.
\bibitem{} Sulewski, P. (2020). Modified Lilliefors goodness-of-fit test for normality. {\em Communications in Statistics-Simulation and Computation,} 3, 1199--1219.

    \bibitem{} Thode, H. C. (2002). \textit{Testing for Normality}, Marcel Dekker, New York.

\bibitem{}  Torabi, H., Montazeri, N. H. and A. Grané (2016). A test of normality
based on the empirical distribution function. {\em SORT}, 40, 55--88.



\end{thebibliography}
\end{document}